\documentclass{amsart}

\usepackage{amsmath,amssymb,latexsym,graphicx}

\numberwithin{equation}{section}
\def\pf{{\bf Proof }}
\newtheorem{theorem}{Theorem}[section]
\newtheorem{lemma}[theorem]{Lemma}

\newtheorem{remark}[theorem]{Remark}

\newtheorem{conjecture}[theorem]{Conjecture}

\def\pf{{\bf Proof }}
\usepackage{enumerate}

\begin{document}

\title{DOMAINS WITH RADICAL-POLYNOMIAL $X$-RAY TRANSFORM}



\author[M. Agranovsky]{MARK AGRANOVSKY}
\address{Bar-Ilan University and Holon Institute of Technology}
\email{agranovs@math.biu.ac.il}





\maketitle

\begin{abstract}

Let $K$ be a compact convex body in $\mathbb R^n.$ For any affine line $L,$ denote $\widehat{\chi}_K(L)=\int_{L}\chi_K(x)dl(x),$ where $dl$ is the arc length measure, the $X$-ray transform of the characteristic function $\chi_K,$ i.e., the length of the chord $K \cap L.$ We prove that if $K$ is bounded by a $C^{\infty}$
real algebraic hypersurface $\partial K$ and the $X$-ray transform $\widehat{\chi}_K(L)$ behaves, under small parallel translations of the line $L$ to the distance $t,$ as the $m$-th root  of a polynomial of $t$, for some fixed $m \in \mathbb N,$ then $\partial K$ is an ellipsoid.

\end{abstract}


\section{Introduction}\label{S:Intro}

This article is devoted to characterization bodies in $\mathbb R^n$ in integral-geometric terms and is motivated by study of so called polynomially integrable bodies. Let us explain this relation.

Given a bounded domain $K \subset \mathbb R^n,$ denote $A_K(\xi,t), \xi \in \mathbb R^n, |\xi|=1, t \in \mathbb R,$ {\it the sectional volume function}, which equals to the $(n-1)$-dimensional volume of the cross-section of $K$ by the affine hyperplane $\{\langle \xi, x \rangle =t\}.$ Here $\langle \ , \ \rangle$ is the inner product in $\mathbb R^n.$

In other words, $A_K(\xi,t)$ is the {\it Radon transform} of the characteristic function $\chi_K:$
$$A_K(\xi,t)=\int\limits_{ \langle \xi, x \rangle=t } \chi_K(x) dvol_{n-1}(x)=vol_{n-1} \big( K \cap \{\langle \xi, x \rangle =t \}\big) .$$
The body $K$ is called {\it polynomially integrable} \cite{Ag}, if $A_K(\xi,t)$ is a polynomial with respect to $t$ so long as the above cross-section is non-empty.

It was proved in \cite{KMY} (see also \cite{Ag}, \cite{Ag1}) that the only polynomially integrable domains with $C^{\infty}$ boundary are solid ellipsoids in odd-dimensional spaces. In particular, the sectional volume function $A_K(\xi,t)$ is never a polynomial in $t$ when $n$ is even.
Nevertheless, for any ellipsoidal domain $E \subset \mathbb R^{2k},$ the {\it squared} sectional volume function $A_E^2(\xi,t)$ does polynomially depend on $t.$  Thus, in any dimension, the sectional volume functions of ellipsoids are either polynomials or radicals of polynomials with respect to $t.$ We conjecture that this property fully characterizes ellipsoids, disregarding the parity of the dimension of the space:

\begin{conjecture}  Let $K \subset \mathbb R^n$ be a compact body with $C^{\infty}$ boundary $\partial K.$ Suppose that
for some $m \in \mathbb N$ the $m$-th power $A_K^m(\xi,t)$ of the sectional volume function is a polynomial in $t$, whenever $A_K(\xi,t) \neq 0.$ Then $\partial K$ is an ellipsoid.
\end{conjecture}
If $\partial K$ is an ellipsoid, then, in the case $n$ is odd, the function $A_K(\xi,t)$ is a polynomial with respect to $t$ of degree $\frac{n-1}{2},$  i.e.,  the condition is satisfied with  $m=1,$ while if $n$ is even, then $A_K(\xi,t)$ is the square root of a polynomial of degree $n-1$ and therefore the condition is fulfilled with $m=2.$
\begin{remark} \label{R:R}
If Conjecture is true then the similar version using $k$-dimensional affine cross-sections, $1 \leq k \leq n-1$ is fixed, is true. It
immediately follows by applying Conjecture to intersections of $K$ with $k+1$-dimensional affine hyperplanes.
\end{remark}

In Theorem \ref{T:Main2} of this article, we  confirm Conjecture for $n=2$ and under  {\it a priori} assumption of algebraicity of the boundary $\partial K.$ Applying this result to two-dimensional sections yields a characterization of $n$-dimensional ellipsoids in terms of the chord length function $\widehat{\chi}_K(L+t\xi),$ i.e., to the situation corresponding in Remark \ref{R:R} to arbitrary $n$ and $k=1.$

Let us start with discussion of the basic, two-dimensional, case.
When $n=2$ then the hyperplanes are affine straight lines $L_{\xi,t}=\mathbb R \cdot \xi^{\perp}+t\xi=
\{x \in \mathbb R^2: \langle x, \xi \rangle =t\} $ and the sectional volume function $A_K(\xi,t)$ boils down to the {\it chord length function}
\begin{equation}\label{E:chord}
A_K(\xi,t) = \widehat{\chi}_K ( L_{\xi,t})=\mbox{{\it length of the chord}} \ K \cap L_{\xi,t}.
\end{equation}

We want to characterize those domains $K$ for which $\widehat \chi_K(\xi,t)$ is an algebraic function  of a simple form, namely, is a radical  of a polynomial in $t.$

There is a relation of the question under discussion with a well known Newton's Lemma about ovals (see \cite{ArVas}). It says that the area cut off a planar domain $K$ with smooth boundary by a straight line is never algebraic function of the parameters of the secant line. The area
$V_K^{\pm}(\xi,t)= \mbox{{\it area}} \ \big (K \cap  \{ \langle x, \xi \rangle \gtrless t \}  \big) $ of a portion of $K$ on one side of the line ({\it the solid volume function}) is just the primitive function of the chord length function. Therefore,
if $A_K(\xi,t)=\sqrt[m]{P_{\xi}(t)}$ where $P_{\xi}(t)$ is a polynomial in $t$ then
\begin{equation}\label{E:Abel}
V_K(\xi,t)=\int \sqrt[m]{P_{\xi}(t)} dt
\end{equation}
is an Abelian integral.

Thus, Newton's lemma says that the solid area function $V_K$ of a domain $K$ in the plane is always transcendental, Theorem \ref{T:Main2} specifies that among those transcendental functions, Abelian integrals (\ref{E:Abel}) characterize ellipses.

Multi-dimensional generalization of Newton's Lemma are related to Arnold conjecture about algebraically integrable domains in $\mathbb R^n$
 (see \cite{Ar},\cite{ArVas},\cite{Vas}).

\section{Main result}

Let $K$ be a compact (connected) domain in $\mathbb R^n.$ Given an affine line $L \subset \mathbb R^n$ and a unit vector $\xi \in \mathbb R^n,$ we will call the function
$$\mu_{L,\xi}(t)=\widehat{\chi_K}(L+t\xi)=\mbox{{\it length of}} \  K \cap (L+t\xi)  , \ \  t \in \mathbb R,$$
{\it the chord length function}.

We will be considering domains $K$ which boundary $\Gamma=\partial K$ is a {\it semi-algebraic} curve. This means that $\Gamma$ is a connected component
of the zero locus of a polynomial $Q$ with real coefficients. The polynomial $Q$ is assumed irreducible over the field $\mathbb C.$

\begin{theorem} \label{T:Main2} Let $K \subset \mathbb R^n, \ n \geq 2, $ be a compact convex domains with $C^{\infty}$ semi-algebraic boundary $\partial K.$ Suppose for any fixed affine line $L, \ L \cap int K \neq \emptyset,$ and unit vector $\xi,$  the chord length function has for small $t$  the form $$\mu_{L,\xi}(t)=\sqrt[m]{P_{L,\xi}(t)},$$
where $m \in \mathbb N$ and $P_{L, \xi}(t)$ a polynomial with respect to $t:$
$$P_{L,\xi}(t)=\sum\limits_{j=1}^N a_j(L, \xi)t^j.$$
Then $\partial K$ is an ellipsoid (and therefore  a posteriori $m$ can be taken 2 ).
\end{theorem}
\begin{remark} The convexity of the domains in Theorem \ref{T:Main2} can be derived from the main condition
 for the chord length function (see \cite{Ag}) and hence Theorem \ref{T:Main2} is valid without assumption of convexity. However, for the sake of simplicity of the exposition, we a priori assume the domain $K$ to be convex.
\end{remark}
Theorem \ref{T:Main2} is, in a certain sense, similar to  Theorem 2 from \cite{IP} which states that if there exists a function $f$ with the X-ray transform identically one then the domain is a ball.

\section{Outline of the proof of Theorem \ref{T:Main2}}
The idea of the proof  is as follows.

First of all, it suffices to prove Theorem \ref{T:Main2} with $n=2,$  then the statement for arbitrary $n$ follows by considering two-dimensional affine cross-sections. In the case $n=2,$ the chord function $\mu_{L,\xi}(t)$ turns to
$\mu_{L,\xi}(t)=A_K(\xi,t)$ if we take $L=\xi^{\perp}=\{x \in \mathbb R^2: \langle \xi, x \rangle=0 \}.$ Also,we will use notation $P(\xi,t)$ instead of $P_{\xi}(t).$

The key point is to determine the degree of the polynomial $t \to P(\xi,t).$ First of all, we want to obtain an upper bound for the degree.
For this purpose it suffices to understand the order of growth of $P(\xi,t)$ as $ t \to \infty.$ Since the information about values of the polynomial $P(\xi,t)$ for large real $t,$ when the line $\{\langle \xi, x \rangle =t\}$ becomes disjoint from  $K,$ is unavailable, we extend $P(\xi,t)$ for complex $t.$

At this point, we use algebraicity of the boundary curve $\partial K$. This curve has a natural complexification, which is a complex algebraic curve in $\mathbb C^2.$ This allows us to construct, in Section 4, analytic extension of the chord length function $A_K(\xi,t)$ to complex values $t \in \mathbb C.$  Determining the growth of the analytic extension along regular pathes going to $\infty$ delivers the upper bound $deg_t P \leq m$ (Section 5).

The lower bound for $deg_t P(\xi,t)$ (Section 6) follows much easier,  from vanishing the chord length function $A_K(\xi,t)$ on tangent lines to $\partial K.$ We show that at Morse points the order of vanishing is $\frac{1}{2}$ and hence $P(\xi,t)=A_K^m(\xi,t)$ vanishes at tangent lines to the order $\frac{m}{2}$. Since there are two tangent lines with the same normal vector $\xi,$ we conclude that  $deg_t P \geq m.$

In Section 7, we finish the proof of Theorem \ref{T:Main2}. Together with the upper bound, this implies $deg_t P=m$ and hence all zeros of $P(\xi,t)$ are delivered by tangent lines. Knowing zeros allows us to reconstruct the polynomial $P(\xi,t)$ up to a factor depending on $\xi,$ and express the chord length function $A_K(\xi,t)$  via the supporting function of $K.$ Then the range conditions (the first three power moments)  for $X$-ray transform applied to the function $A_K$ imply that the supporting function of $K$ coincides with the supporting function of an ellipse.

\section{Analytic continuation of the chord length function}

Let $n=2.$ Let $K$ be a domain satisfying the conditions of Theorem \ref{T:Main2}.
 We assume that the boundary $\Gamma=\partial K$ is a non-singular real semi-algebraic curve, which means that there is a real irreducible polynomial $Q(x_1,x_2)$ such that
$$Q(x)=0, \ x=(x_1,x_2) \in \Gamma$$
and $\nabla Q(x) \neq 0, x \in \Gamma.$

\begin{figure}[h]
\centering
\includegraphics[width=12cm]{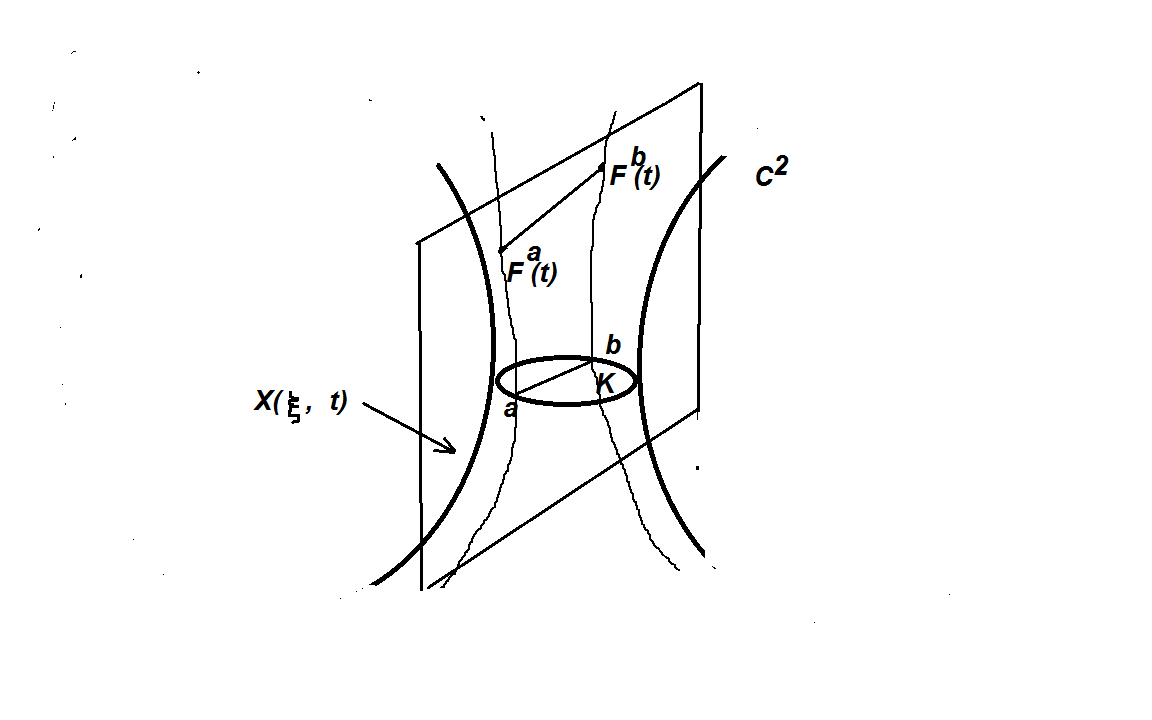}
\end{figure}

Extend polynomial $Q$ to the complex space $\mathbb C^2$ and denote $\Gamma^{\mathbb C}$ the complex algebraic curve
$$\Gamma^{\mathbb C}=\{ z=(z_1,z_2) \in \mathbb C^2: Q(z)=0\}.$$
The domain $K$ is regarded as a set in the real subspace $\mathbb R^2=\{z \in \mathbb C^2: Im z_1=Imz_2=0\},$ so that $\partial K \subset \Gamma^{\mathbb C} \cap \mathbb R^2.$

Given a real unit vector $\xi \in \mathbb R^2$ and $t \in \mathbb C$ denote the complex affine line in $\mathbb C^2:$
$$X(\xi,t)=\{ z \in \mathbb C^2: \langle \xi, z \rangle=\xi_1z_1+\xi_2z_2=t\}.$$

\begin{lemma} \label{L:path}
Fix a unit vector $\xi_0 \in S^1.$  There is a finite set $Z_0 \subset \mathbb R$ such that is $t_0 \in \mathbb R \setminus Z_0$ and the real affine line $\{ \langle \xi_0, x \rangle  =t_0 \}$ intersects transversally the curve $\Gamma=\partial K$ at two points $a$ and $b$ the following is true.
There is a path $T \subset \mathbb C,$ joining $t_0$ and $\infty,$ an open connected neighborhood $U \subset \mathbb C$ of $T$
and two holomorphic mappings
$$F^a, \ F^b : U \to \mathbb C^2$$
such that
\begin{enumerate}[(i)]
\item $F^a(t), \ F^b(t) \in X(\xi_0,t) \cap \Gamma^{\mathbb C}$ for all $t \in U.$
\item If $t \in U \cap \mathbb R$ then $F^a(t), \ F^b(t) \in \Gamma.$
\item $F^a(t_0)=a, \ F^b(t_0)=b.$
\end{enumerate}
\end{lemma}

\pf Applying a suitable rotation in $\mathbb R^2,$ we can assume for simplicity that
$$\xi_0=(0,1).$$
In this case, $X(\xi,t)$ is the complex line $\{z_2 = t\}$ in $\mathbb C^2$ and
$a=(a_1,t), \ b =(b_1,t).$
The condition $ z \in X(\xi_0, t) \cap \Gamma^{\mathbb C}$ translates in this case as $Q(z_1,t) = 0.$

Then Lemma \ref{L:path} asserts, in fact, that if we consider the projection
$$\pi:\Gamma^{\mathbb C} \to \mathbb C, \ \pi(z_1, t)=t$$
then there is a path $T \subset \mathbb C$, joining $t_0$ and $\infty$ and a neighborhood $U$ of $T$ such that the holomorphic mapping $\pi$
possesses two holomorphic sections $F^a, \ F^b: U \to \Gamma^{\mathbb C}$ of $\pi$  over the set $U$ with the initial conditions {\it (ii), (iii)}. The existence of such sections follows from the path lifting property (see, e.g., \cite{Fu}, Proposition 11.6) of covering maps and from the fact that the projection $\pi$ is a covering outside of the finite set of poles and ramification points.

Let us give more extended analytic arguments, for the sake of  self-sufficiency.
Represent the polynomial $Q,$ defining the complex algebraic curve $\Gamma^{\mathbb C},$ in the form
$$Q(z_1,z_2)=q_0(z_2)+q_1(z_2)z_1+ \cdots + q_{M}(z_2)z_1^{M},$$
where $q_j$ are polynomials of one complex variable and $q_M \not\equiv 0.$

Consider the discriminant $D(t)=Disc_{z_1}Q_t$ of the polynomial $Q_t(z_1):=Q(z_1,t):$
$$D(t)= q_M^2 (t) \prod\limits_{i<j}(r_i -r_j)^2,$$
where $r_i$ are the roots of the polynomial $Q_t.$

The discriminant $D(t)$ is a polynomial in coefficients $q_0(t),\cdots, q_M(t)$ and therefore is itself a polynomial in $t.$ It cannot vanish identically. Indeed, if $D(t) \equiv 0$ then for any $t \in \mathbb C$ either $q_M(t)=0$ or the polynomial $Q_t$ has at least one multiple zero, which is a common zero of $Q_t$ and $\frac{\partial Q_t}{\partial z_1}.$ Since the number of zeros of $q_M$ is finite, the set of common zeros of
$Q$ and $\frac{\partial Q}{\partial z_1} $ must be then an infinite set in $\mathbb C^2$. By Bezout theorem this means that the two polynomials have a common polynomial factor, which is impossible as $Q$ is irreducible and $\frac{\partial Q}{\partial z_1}$ is of less degree than $G.$
Therefore $D(t)$ is a nontrivial polynomial and hence the discriminant set
$$Z:=\{t \in \mathbb C: D(t)=0\}$$ is finite.
We set $Z_0 =Z \cap \mathbb R.$

Let $t_0$ is real, $ t_0 \notin Z_0.$ Let the line $\{ \langle \xi_0, x \rangle =t_0 \}$ intersects $\Gamma$
at the points $a$ and $b.$

Now choose a smooth simply-connected path $T \subset \mathbb C \setminus Z,$ joining $t_0$ and $\infty$ and
avoiding the discriminant set $Z.$ Consider the equation
$$Q_t(z_1)= Q(z_1, t)=0, \ t \in T.$$
Since $q_M(t) \neq 0$ for $t \in T,$ we can divide by $q_M(t)$ and reduce the polynomial in the equation to the form:
\begin{equation}\label{E:eq}
Q_t(z_1)= p_0(t)+ p_1(t)z_1 + \cdots + p_{M-1}(t)z_{M-1}+ z_1^M=0,
\end{equation}
where the coefficients
\begin{equation}\label{E:coeff}
p_j(t)=\frac{q_j(t)}{q_M(t)}, \ j=0,\cdots M, \ p_M(t)=1,
\end{equation}
are continuous functions on $T.$

Thus, we deal with an algebraic monic equation for $z_1$ with the coefficients, continuously depending on the parameter $t \in T.$
For any $t \in T$ we have $t \notin  Z,$ hence $D(t) \neq 0,$ i.e., all the roots of $Q_t$ are simple.

The monodromy theorem  \cite{Hu}, Thm. 16.2, implies that the algebraic equation (\ref{E:eq}) so completely solvable. This means that there is no monodromy
on $t \in T$ and there exist $M+1$ continuous functions $f_0(t), \cdots, f_M(t)$ on $T,$ satisfying equation (\ref{E:eq}):
$$Q_t(f_j(t))=0, \ j=0,...,M.$$
Since all the roots are simple, we have $f_i(t) \neq f_j(t)$ for all $t \in T$ and $ i \neq j.$

Let us explain  this point in more details. Consider the spaces
$$E=\{ (\lambda_0, \cdots \lambda_M)
 \in \mathbb C^{M+1}: \ \lambda_i \neq \lambda_j, \ i,j=0,...,M, \}$$
and
$$B=\{ (p_0,...,p_{M-1}) \in \mathbb C^M: p_0+ p_1z_1+...+p_{M-1}z_1^{M-1}+z_1^M \ \mbox{{\it has no multiple roots}} \}.$$
Define the mapping
$$ \pi: E \to B$$
as follows: $\pi(\lambda) , \lambda \in E$ is the vector $p=(p_0,...,p_{M-1})$ of the coefficients of the monic polynomial with the roots $\lambda_j,$ i.e.,
$$ \sum\limits_{j=0}^{M}p_j \lambda^j =\prod\limits_{j=0}^{M}(\lambda-\lambda_j), \ p_M=1.$$
The roots $\lambda_j$ of the polynomial in the right hand side are symmetric functions of the coefficients $p_j$  and
by Implicit Function Theorem,  $\pi$ is a $(M+1)!$-covering map.

Let $\lambda^{(0)}=(\lambda_0^{(0)}, \cdots, \lambda_M^{(0)}) \in E$ be the roots of the polynomial $Q_{t_0}(z_1)=Q(z_1,t_0),$ i.e.,
$\pi(\lambda^{(0)})= \big(p_0(t_0), ..., p_{M-1}(t_0) \big).$

The mapping
$$g(t):= \big(p_0(t),...,p_{M-1}(t)\big) \in B, \ t \in T,$$
where $p_j(t)$ are coefficients (\ref{E:coeff}) of the polynomial $Q_t$
defines a path
$$g: T \to B$$ in the base space $B.$ The path lifting property of covering mappings (see, e.g., \cite{Fu},Proposition 11.6; \cite{Hu},Theorem 16.2) says that there is a lifting path
$$f=(f_0, \cdots, f_{M}): T \to E$$ such that $\pi \circ f= g$ and $f(t_0)=\lambda^{(0)}.$

Then the functions
$$f_j(t), j =0, \cdots, M$$
define the above claimed continuous family of roots of the polynomials $Q_t.$

The point $a=(a_1,t_0), \ b=(b_1,t_0) $ satisfy $Q(a_1,t)= Q(b_1,t)=0.$ Therefore, there are two branches, say, $f_j(t), \ f_j(t)$ which take at $t_0$ the values $a_1, \ b_1,$ correspondingly. Denote $f_i=f^a, \ f_j=f^b.$

Then we have
\begin{enumerate}
\item $Q(f^a(t),t)= Q_t(f^a(t)) = 0, \ t \in T,$
\item $Q(f^b(t),t)= Q_t(f^b(t)) = 0, \ t \in T,$
\item $f^a(t_0)=a_1, \ f^b(t_0) = b_1.$
\end{enumerate}

For any fixed $t \in \mathbb C$ we have
$$Q(f^a(t),t)=0, \ \frac{\partial Q}{\partial z_1}(f^a(t),t) \neq 0,$$
because $f^a(t)$ is a simple root. The same is true for $f^b(t).$
By Implicit Function Theorem, there is a complex neighborhood $U^a_t$ of $t$ and a complex  neighborhood $V_{f^a(t)}$ of $f_a(t)$ such that
for any $s \in U^a_t$ there exists a unique $w:=f^a(s) \in V_{f^a(t)}$ such that $Q(w,s)=0,$ and the function $w=f^a(s)$ is holomorphic in $U^a_t.$

Thus, given $t \in T$ the function $f^a(t)$ extends to the neighborhood $U^a_t$ as a holomorphic function. The union $U^a=\cup_{t \in T} U^{a_t}$ constitutes an open connected set containing the path $T.$ The function $f^a(t), t \in T$ extends to $U^a$ as a holomorphic function.
The extensions satisfies the same polynomial equation $Q(f^a(t),t)=0, \ t \in U^a.$  Similarly, we construct an open set $U^b$
and the holomorphic function  $f^b$ in $U^b$ with analogous properties.

Now set
\begin{equation}\label{E:U-F}
\begin{aligned}
&U=U^a \cap U^b,\\
&F^a(t)=(f^a(t),t), \ F^b(t)=(f^b(t),t).
\end{aligned}
\end{equation}

Check properties {\it (i)-(iii).}
The mapping $F^a: U \to \mathbb C^2$ is holomorphic in $U,$  and by construction $Q(F^a(t))= Q(f_a(t),t) = 0$ for all $t \in U.$  Also, since $\xi=(0,1),$ we have $\langle \xi, F^a(t) \rangle = t,$ hence $F^a(t) \in \Gamma^{\mathbb C} \cap X(\xi,t).$

Furthermore, $F^a(t_0) = (f^a(t_0), t_0) = (a_1, t_0) = a.$  For  $t \in U \cap \mathbb R$ near $t_0,$   the straight line $x_2=t$ intersects $\partial K$ at a point $a_t = (a_{1,t},t)$ close to $a.$ Then $Q(a_{1,t},t)=0$ for $t$ in a neighborhood of $t_0,$ i.e. the polynomial $Q$ vanishes on an open subarc of the real-analytic curve $\partial K.$ By the uniqueness theorem for holomorphic functions, the identity holds for all real  $t \in U.$  Since the root is unique, we conclude that $a_t=F^a(t)$ and thus $F^a(t) \in \partial K \cap U. $
Hence the set $U$ and the mapping $F^a$ satisfy all properties {\it (i)-(iii)}. Similarly, we proceed with the mapping $F^b.$

Lemma is proved.

\begin{lemma}\label{L:ext} Assume, as in Theorem \ref{T:Main2}, $A^m_K(\xi,t)=P(\xi,t),$ where $P(\xi,t)$ is a polynomial in $t.$
Fix a unit vector $\xi_0 \in S^1.$ Let the straight line $\langle \xi_0, x \rangle =t_0$ intersects $\partial K$  at the points $a$ and $b$ and $t_0 \notin Z_0,$ where $Z_0$ is the finite exceptional set from Lemma \ref{L:path}.
Construct the open set $U \subset \mathbb C$ and the holomorphic mappings $F^a, F^b: U  \to \mathbb C^2$ as in Lemma \ref{L:path}.
Then
\begin{equation}\label{L:P=F}
P(\xi_0,t)= \big (\langle \xi_0^{\perp}, F^a(t)-F^b(t) \rangle \big)^m, \ t \in U.
\end{equation}
where $\xi_0^{\perp}=\frac{a-b}{|a-b|}$ is the vector orthogonal to $\xi_0.$
\end{lemma}

\pf By Lemma \ref{L:path}, {\it (i), (ii)}, when $t \in U \cap \mathbb R$ then the segment $[F^a(t), F^b(t)]$ is just the chord $X(\xi_0, t) \cap K.$
The length of the chord is
$$A_K(\xi,t)=|F^a(t)-F^b(t)|=\langle F^a(t)-F^b(t), \frac{F^a(t)-F^b(t) }{|F^a(t)-F^b(t)|} \rangle.$$
Also, $\langle F^a(t)-F^b(t), \xi_0 \rangle = t-t=0$ and hence the second factor in the inner product is a unit vector orthogonal to $\xi_0$ and does not depend on $t.$ By taking $t=t_0,$ we find, due to $F^a(t_0)=a, \ F^b(t_0)=b,$ that this vector is $\xi_0^{\perp}=\frac{a-b}{|a-b|}.$

Then $P(\xi_0,t)=A^m(\xi_0,t)= \big (\langle F^a(t)-F^b(t), \xi_0^{\perp}  \rangle \big)^m $ for $t \in U \cap \mathbb R.$
Since $P(\xi,t)$ and the function in the right hand side of the equality  are holomorphic in $t \in U,$ and the set $U$ is open and connected, the equality is satisfied for all $ t \in U$ by the uniqueness theorem. Lemma is proved.

\section{Upper bound of the degree of the polynomial $P(\xi,t)$}

As before, $n=2.$ We fix $\xi_0 \in S^1$ and  $t_0 \in \mathbb R$ such that the line $\{\langle \xi_0, x \rangle $ meets $\Gamma=\partial K$  at the points $a,b$ and $t$ does not belong to the finite exceptional set $Z_0$ in Lemma \ref{L:path}. Let the path $T \subset \mathbb C$ and the open set $U \subset \mathbb C, \ T \subset U,$ be as in Lemma \ref{L:path}.

In order to obtain  the upper bound for the degree of the polynomial $P(\xi_0, t),$ it suffices to estimate the order of its growth as $t \to \infty.$ According to the representation (\ref{L:P=F}), the problem is reduced to understanding the behaviour of the mappings $F^a(t), F^b(t)$ as $t \to \infty.$

Let
$$Q=Q_0+\cdots +Q_N$$
be the decomposition of the polynomial $Q,$ defining the complexified boundary $\Gamma^{\mathbb C},$
into homogeneous polynomials $Q_j, deg Q_j=j, j=1,..., N=deg Q.$

The leading homogeneous polynomial $Q_N(z_1,z_2)$ of two complex variables is completely reducible over $\mathbb C:$
$$Q(z_1,z_2)=const \prod\limits_{j=1}^J (A_j z_1+B_jz_2)^{m_j}, \ \sum\limits_{j=1}^J m_j=N=deg Q.$$

\begin{lemma} \label{L:growth} Suppose that $\alpha_0 B_j - \beta_0 A_j \neq 0$ for any $j=1,\cdots,J,$
where $\xi_0=(\alpha_0, \beta_0).$
Then $\frac{F^a(t)}{t}, \ \frac{F^b(t)}{t}$ are bounded for $t \in U.$
\end{lemma}

\pf The zero set $Q_N^{-1}(0)$ of the homogeneous polynomial $Q_N$ in $\mathbb C^2$ consist of $J$ complex lines $A_j z_1+B_j z_2=0, j =1,...,J,$ counting multiplicities.
The condition for $\xi_0$ means that the affine complex line $X(\xi_0,0)$ is none of those.
Again, it would be convenient to apply rotation in the real plane $\mathbb R^2$ and make
$\xi_0=(0,1).$ Then we have $\alpha_0=0, \beta_0=1$ and $A_j \neq 0$ for all $j=1,...,J.$

Let $z \in Q^{-1}(0) \cap X(\xi_0,t).$ It means that $z_2=t$ and $Q(z)=Q(z_1,t)=0,$  or, the same,
$$Q_0 + tQ_1(\frac{z_1}{t}, 1)+...+ t^{N-1} Q_{N-1}(\frac{z_1}{t}, 1)+t^N Q_N(\frac{z_1}{t}, 1)=0.$$
Dividing both sides by $t^N$ yields
\begin{equation}\label{E:equality}
\Psi(w,s):=s^N  Q_0+ s^{N-1} Q_1(w,1)+ \cdots + s Q_{N-1}(w,1)+ Q_N(w,1) =0,
\end{equation}
where $w, s$ and $z_1, t$ are related by
$$w=\frac{z_1}{t}, \ s=\frac{1}{t}.$$
Denote
$$w_j=-\frac{B_j}{A_j}$$
the root of the polynomial $Q_N(w,1),$ of multiplicity $m_j.$

Consider the logarithmic residue
$$r_j(s)=\frac{1}{2\pi i}\int\limits_{|w-w_j|=\varepsilon} \frac {\Psi_{w}^{\prime}(w,s)}{\Psi(w,s)} dw,$$
where $\varepsilon$ is so small that $w_j$ is the only root of $\Psi(w,0)=Q_N(w,1)$ in the disc $|w-w_j| \leq \varepsilon.$
Then
$$r_j(0)=m_j.$$
Thus, the total sum is
$$r_1(0)+ \cdots + r_J(0)=m_1+ \cdots +m_J=N.$$

By continuity, for every $j=1,...,J$ and sufficiently small $\varepsilon >0$ there exists $\delta_j >0$ such that for $|s|<\delta_j,$  $\Psi(w,s) \neq 0$ when $|w-w_j|=\varepsilon.$
Then the function $r_j(s)$ is continuous in $|s| <\delta_j$ and integer-valued, therefore $r_1(s)+\cdots +r_J(s)=N$ for
$|s|<\delta=\min \{\delta_j,...,\delta_J \}.$ . This means that for $|s|<\delta,$ all $N$ roots, counting multiplicities, of the polynomial $w \to \Psi(w,s)$  are located in $\varepsilon$-neighborhood of the set of zeros  of the polynomial $\Psi(w,0)=Q_N(w,1).$

Going back to the variable $z_1=tw$  and the function $Q(z_1, t)$
we conclude that if $|t|>\frac{1}{|\delta|},$ then $|\frac{z_1}{t}-w_j| < \varepsilon$ for some root $w_j$ of $Q_N(w,1).$

Since we consider the situation $\xi_0=(0,1),$  the condition  {\it (i)} in Lemma \ref{L:path}  means that
$F^a(t), F^b(t)$ have the form
$$F^a(t)=\big ( f^a(t), t \big), \ F^b(t)=\big (f^b(t), t \big), t \in U.$$

Now, by the construction in Lemma \ref{L:path}, if $t \in U$ then $Q_t(z_1)=Q(z_1, t)$ has only simple roots
$\Lambda(t)=\{\lambda_0(t),...,\lambda_M(t)\}.$ Therefore, when $|t|> \frac{1}{\delta}$ then $\frac{\lambda_j(t)}{t}$ are in an $\varepsilon$-neighborhood of a root of $Q_N(w,1).$  Among the collection  $\Lambda(t)$ of the roots of $Q_t$, there are two branches, say, $\lambda_i(t), \lambda_j(t),$ determined by the initial conditions $\lambda_i(t_0)=a, \lambda_j(t_0)=b.$  They must coincide, correspondingly, with $f^a(t)$ and $f^b(t)$ and hence
$$  \bigg| \frac{f^a(t)}{t}-w_i \bigg| <\varepsilon, \ \bigg| \frac{f^b(t)}{t}-w_j \bigg| < \varepsilon, \ t \in U,  |t| > \frac{1}{\delta}.$$
Since the set of roots $w_i$ is finite, we conclude that $\frac{f^a(t)}{t}, \frac{f^b(t)}{t}, t \in U,$ are bounded. Lemma is proved.

\begin{lemma}\label{L:above} For all $\xi \in S^1$ but finite set of those,
$\deg_t P(\xi,t)$ is at most $m.$
\end{lemma}
\pf The set of vectors $\xi_0$ which do not satisfy the condition of Lemma \ref{L:growth} is finite.
Let $\xi=\xi_0$ be not such a vector and let $t_0, a, b,$ be as in Lemma \ref{L:path}.
By Lemma \ref{L:ext}, $A^m_K(\xi,t)=P(\xi,t)=\big( \langle \xi^{\perp}, F^a(t)-F^b(t) \rangle \big)^m.$
Then by Lemma \ref{L:growth} $\frac{P(\xi,t)}{t^m}$ is bounded, as $ t \to \infty, \ t \in U.$
Therefore, $deg_t P(\xi,t) \leq m.$

\section{The zeros and exact degree of $P(\xi,t)$}

In this section $n=2.$ Given $\xi \in S^1,$ denote
$$  \rho_{+}(\xi)=\max\limits_{x \in K} \langle \xi, x \rangle, \ \rho_{-}(\xi)=\min\limits_{x \in K} \langle \xi, x \rangle.$$
Denote $M_{\pm}(\xi) \in \partial K$ the points where the maximum and minimum are attained:
$$\rho_+(\xi)=\langle \xi, M_+(\xi) \rangle,  \  \rho_{-}(\xi)=\langle \xi, M_{-}(\xi) \rangle.$$
The lines $\langle \xi, x \rangle =\rho_{\pm}(\xi)$ are supporting lines to $\partial K$ at the points $M_{\pm}(\xi)$ and the vectors $\pm \xi$ are the unit outward normal vectors $\pm \xi =\nu_{M_{\pm}(\xi)}$ to the curve $\partial K$ at the points $M_{\pm}(\xi),$ correspondingly.

\begin{lemma}\label{L:below}
The polynomial $P$ has the form
$$P(\xi,t)=c(\xi)\big(\rho_{+}(\xi)- t \big)^{\frac{m}{2}} \big(t-\rho_{-}(\xi) \big)^{\frac{m}{2}}, \ c(\xi)>0.$$
\end{lemma}
\pf By continuity, it suffices to establish the representation for almost all $\xi.$ Since $\partial K$ is a non-singular algebraic curve, it is real-analytic. Then all points, except finite number of those, are Morse points. Therefore, only for a finite vectors $\xi,$ the points $M_{\pm}(\xi)$ are non-Morse. Choose $\xi$ such that this is the case.

Using rotation and translation, we can assume that $\xi=(0,1)$ and $M_(\xi)=(0,0).$ The outward normal vector at ${M_{-}(\xi)}$ is $-\xi=(0,-1)$ and
$\rho_{-}(\xi)=0.$
In a neighborhood of the point $M_{-}(\xi)=(0,0)$ the curve $\Gamma$ can be represented as
a graph $x_2=\varphi(x_1)$ of a real-analytic function, with $\varphi^{\prime}(0,0)=0, \varphi^{\prime\prime}(0,0) \neq 0:$
$$x_2=\frac{1}{2} \varphi^{\prime\prime}(0,0)x_1^2+ o(x_1), \ x_1 \to 0.$$
Then the length of the chord $x_2=t$ is $A_K(\xi,t)= 2\sqrt{\frac{2t}{\varphi^{\prime\prime} (0,0) } } + o (\sqrt{t}), t \to 0.$

It shows that the length of the chord, obtained by the parallel translation of a tangent line to the distance $t$,
behaves at $\sqrt{t}.$  This yields that if  $M_{\pm}(\xi) \in \Gamma$ are Morse points then
$$A_K(\xi,t)=const (t- \rho_{-}(\xi))^{\frac{1}{2}}+o( (t -\rho_{-}(\xi))^{\frac{1}{2}}, t \to \rho_{-}(\xi).$$
Similarly, $$A_K(\xi,t)=const (\rho_{+}(\xi)-t)^{\frac{1}{2}}+o( (\rho_{+}(\xi) -t)^{\frac{1}{2}}, t \to \rho_{+}(\xi).$$
Thus, the polynomial $P(\xi,t)=A^m_K(\xi,t)$ vanishes, to  the order $\frac{m}{2},$
at  $t=\rho_+(\xi)$ and $t=\rho_{-}(\xi).$   This means, firstly, that $m$ is even and, secondarily, that $deg_t P(\xi,t) \geq m.$
According to the remark at he beginning of the proof, it holds for all $\xi \in S^1$ except for a finite set.
Together with Lemma \ref{L:above} it proves that $deg_t P(\xi,t)=m$ for all $\xi$ except for  a finite set, and for those $\xi$ we have
$P(\xi,t)=c(\xi)\big(\rho_+(\xi)-t \big)^{\frac{m}{2}} \big(t - \rho_{-}(\xi) \big )^{\frac{m}{2}}.$ By continuity, $P(\xi,t)$ has the claimed representation for all $\xi \in S^1.$

\section{Proof of Theorem \ref{T:Main2}}

First of all, as it has been mentioned before, it suffices to prove Theorem \ref{T:Main2} for $n=2.$ Indeed, each transversal intersection of $K$ with two-dimensional affine plane produces a domain in this plane satisfying all the conditions of Theorem \ref{T:Main2}. If we could conclude that all such two-dimensional cross-sections are bounded by ellipses then the entire body $K$ is bounded by an ellipsoid.

Thus, we assume that $n=2.$ Then we follow  the arguments from \cite{Ag}.
The representation of the polynomial $P(\xi,t)$ given by  Lemma \ref{L:below} yields
\begin{equation}\label{E:A}
A_K(\xi,t)= \sqrt[m]{P(\xi,t)} = d(\xi) \big(\rho_{+}(\xi) -t \big)^{\frac{1}{2}} \big (t-\rho_{-}(\xi)\big) ^{\frac{1}{2}} ,
\end{equation}
where $d(\xi)=\sqrt[m]{c(\xi)}.$
Representation (\ref{E:A}) holds whenever $t \in [\rho_{-}(\xi), \rho_+(\xi)],$ otherwise $A_K(\xi,t)=0.$

The next step is applying the range conditions for Radon transform \cite{He}. Function $A_K(\xi,t)$ is the Radon transform of the characteristic function $\chi_K$ and hence satisfy the moment conditions. Namely, the moments
\begin{equation}\label{E:M}
M_k(\xi)=\int\limits_{\rho_{-}(\xi)}^{\rho_+(\xi)} A_K(\xi,t)t^k dt
\end{equation}
must be restriction to the unit circle $S^1$ of a homogeneous polynomial of degree $k.$
Notice, that $M_0(\xi)= \mbox{ {\it area}} \ K =const >0.$

Then substituting (\ref{E:M}) into (\ref{E:A})  yields:
$$M_k(\xi)=d(\xi)\int\limits_{\rho_{-}(\xi)}^{\rho_{+}(\xi)} \big(\rho_{+}(\xi) -t \big)^{\frac{1}{2}} \big (t-\rho_{-}(\xi)\big) ^{\frac{1}{2}}t^k dt.$$
Denote
\begin{equation}\label{E:BC}
\begin{aligned}
&B(\xi)=\frac{1}{2} \big (\rho_{+}(\xi) + \rho_{-}(\xi) \big),\\
&C(\xi)=\frac{1}{2} \big (\rho_+(\xi)-\rho_{-}(\xi) \big).
\end{aligned}
\end{equation}
and perform the change of variables in the integral
$$s=t-B(\xi).$$
Since
$$\big(\rho_{+}(\xi)-t\big) \big(t-\rho_{-}(\xi)\big)=\big(C(\xi)^2 -s^2 \big),$$
substitution in the integral yields:
$$M_k(\xi)=d(\xi) \int\limits_{-C(\xi)}^{C(\xi)}  \big ( C(\xi)^2 -s^2 \big) ^{\frac{1}{2}} \big(s+B(\xi)\big)^kds.$$

Finally, the change of variable
$$s=C(\xi)v$$ leads to
$$M_k(\xi)=G(\xi) \int\limits_{-1}^{1}(1-v^2)^{\frac{1}{2}} \big( C(\xi)v + B(\xi) \big) ^k  dv$$
where
$$G(\xi)=d(\xi)C(\xi).$$
Let us write the first three moments:
\begin{equation}\label{E:moments}
\begin{aligned}
&M_0(\xi)=G(\xi)\alpha_0, \\
&M_1(\xi)=G(\xi) \big( C(\xi)\alpha_1+B(\xi) \alpha_0 \big),\\
&M_2(\xi)=G(\xi) \big(C^2(\xi)\alpha_2 +2C(\xi)B(\xi)\alpha_1+ B^2(\xi)\alpha_0 \big),
\end{aligned}
\end{equation}
where
$\alpha_k=\int\limits_{-1}^1 (1-v^2)^{\frac{1}{2}}v^k dv.$

Since $\alpha_1=0$ we have
$$B(\xi)=\frac{M_1(\xi)}{M_0(\xi)}.$$
By the range conditions, $M_0(\xi)=const, G(\xi)=const$ on $S^1$  and  $M_1(\xi)$ extends from $S^1$ as a homogeneous linear polynomial.
Therefore, $B(\xi)$ is the restriction to the unit circle  of a linear form:
$$B(\xi)=\langle \xi, b \rangle,$$
where $b$ is a fixed vector.

Since $M_2(\xi)$ extends from $S^1$ as a homogeneous quadratic polynomial, $G(\xi)$ is constant and $B(\xi)$ is a linear form, the third equality in (\ref{E:moments}) implies that $C^2(\xi)$ is the restriction of a (strictly positive) quadratic form.
Applying a suitable rotation we can reduce $C^2(\xi)$ to the form
\begin{equation}\label{E:C}
C^2(\xi)=c_1 \xi_1^2+c_2\xi_2^2, \ c_1,c_2 > 0.
\end{equation}

Now, by translating $K$ by the vector $b,$ we can make $B(\xi)=0$ for all $\xi.$ Indeed, the supporting functions $\rho_{\pm}$ transform, under the translation by the vector $b,$ as follows:
$$\rho_+(\xi) \to \rho_+(\xi)- \langle \xi, b \rangle, \   \rho_{-}(\xi) \to \rho_{-}(\xi)- \langle \xi, b \rangle.$$
Then (\ref{E:BC}) shows that $B(\xi)$ transforms to $B(\xi) - \langle \xi, b \rangle =0.$

Thus, we can apply the translation $K$ by the vector $b$ and assume that $B(\xi)=0.$
Then from (\ref{E:BC}) we have $\rho_+(\xi)=-\rho_{-}(\xi)$ and $C(\xi) =\rho_+(\xi).$
From (\ref{E:C}), we obtain
$$\rho_+(\xi)=\sqrt{c_1 x_1^2+ c_2 x_2^2}.$$
Thus, the supporting function $\rho_{+}(\xi)$ of the body $K$  coincides with the supporting  function of the ellipse $$E=\{ \frac{x_1^2}{c_1}+\frac{x_2^2}{c_2}=1\}.$$ Thus, we conclude that $\partial K=E.$ Theorem \ref{T:Main2} is proved.



\noindent
Bar-Ilan University and Holon Institute of Technology;  Israel.

\noindent
{\it E-mail address:} agranovs@math.biu.ac.il

\end{document}